\documentclass[12pt,twoside]{article}
\usepackage{amssymb,amsmath,mathrsfs,txfonts,graphicx,color}
\usepackage{epstopdf}
\usepackage{cite}
\usepackage{bbm}
\usepackage{bm}
\usepackage{amsfonts}
\usepackage{cite}
\usepackage{multirow}
\usepackage{float} \usepackage[colorlinks=true]{hyperref}
\let\oldsection\section
\renewcommand\section{\setcounter{equation}{0}\oldsection}

\newtheorem{theorem}{Theorem}[section]
\newtheorem{lemma}{Lemma}[section]
\newtheorem{proposition}{Proposition}[section]
\newtheorem{definition}{Definition}[section]
\newtheorem{remark}{Remark}[section]

\allowdisplaybreaks

\begin{document}

\title{\Large\bf Global boundedness and asymptotic stability of the Keller-Segel
system with logistic-type source in the whole space}
\author {Qingchun Li,
Haomeng Chen\thanks{Corresponding author. {\it E-mail}:  15100616217@163.com}
}
\date{}

\maketitle

\begin{abstract}
In this paper, we investigate the Cauchy problem of the parabolic-parabolic Keller-Segel system
with the logistic-type term $au-bu^\gamma$ on $\mathbb{R}^N, N\geq2$. We discuss the global boundedness of classical solutions with
nonnegative bounded and uniformly continuous initial functions when $\gamma>1$.
Moreover, based on the persistence of classical solution we show the large time behavior of
the positive constant equilibria with strictly positive initial function in the case of $\gamma\in(1,2)$.
\end{abstract}


\section{Introduction}
Chemotaxis, the directed movement of cells, bacteria, single cells,
and multicellular organisms in response to specific chemical gradients,
holds significant implications across a diverse range of biological
processes \cite{TH}.
The phenomenon of chemotaxis was initially observed in the 1880s
during a photosynthesis experiment on water cotton by German botanist
Engelmann and microbiologist Pfeffer \cite{TE,WP1,WP2}.
In separate studies, they discovered that certain algae engage in
oxygen production through photosynthesis under light conditions,
while bacteria tend to aggregate at elevated levels of oxygen concentration.
After that, the chemotactic mechanism of bacteria has garnered significant
scientific attention, thereby fostering extensive research on bacterial chemotaxis.

In 1970, Keller and Segel established the first mathematical chemotaxis model
to describe the aggregation phenomenon of the slime mold Dyctyostelium discoideum
\cite{EK1,EK2,EK3}:
\begin{align} \label{1.2}\allowdisplaybreaks
 \left\{
    \begin{array}{llll}
        \displaystyle u_t=\Delta u-\chi\nabla\cdot(u\nabla v), &&x \in\Omega, t>0,
        \\
        \displaystyle v_t=\Delta v-v+u, &&x \in \Omega, t>0,
    \end{array}
 \right.
\end{align}
where $\Omega\subset\mathbb{R}^N$ is a smooth and bounded domain.
Here $u(x,t)$ and $v(x,t)$ denote the density of mobile cells and
the concentration of chemical signal at the position $x$ and time $t$, respectively.
Biologically, the constant $\chi>0$ shows the intensity of chemotaxis.
The model \eqref{1.2} has been extensively investigated over the
past few decades. A large variety of work has been dedicated to
the phenomenon that solutions may exhibit global boundedness
or finite-/infinite-time blow-up, with the occurrence of blow-up
being dependent on the space dimension.
In fact, the occurrence of finite-/infinite-time blow-up of the problem
\eqref{1.2} is not observed for $N=1$.
In this case, Osaki and Yagi \cite{KO1} proved the global boundedness of classical
solutions by establishing a priori estimates for the local solutions.
If $N=2$, the blow up of solutions is determined by the initial mass of
bacteria $\int_\Omega u_0$. For the radially symmetric domain $\Omega$,
Nagai et al. \cite{TN2} showed that the global solution is uniformly bounded
when $\int_\Omega u_0<\frac{8\pi}{\chi}$. In addition, Herrero and Velazquez \cite{MA} proved
that the solution blows up at finite time when $\int_\Omega u_0>\frac{8\pi}{\chi}$.
For other asymmetric domain $\Omega$, if $\int_\Omega u_0<\frac{4\pi}{\chi}$,
the boundedness of global classical solutions was obtained by Nagai et al. \cite{TN2}.
If $\int_\Omega u_0>\frac{4\pi}{\chi}$, Horstmann and Wang \cite{DH} established the blow-up
solution of the system \eqref{1.2}.
In the case of $N\geq3$, Winkler \cite{MW1} showed that the solution of the system
\eqref{1.2} occurs blow-up phenomenon in the finite-time.
Moreover, Winkler \cite{MW2} derived the global existence and boundedness of solutions
for sufficiently small initial data satisfying $\|u_0\|_{L^p(\Omega)}(p>\frac{N}{2})$
and $\|\nabla v_0\|_{L^q(\Omega)}(q>N)$.
Similarly as above, corresponding results have also been obtained in two
and higher dimensional environments when $\Omega=\mathbb{R}^N$.
In the case of $N=2$, Herrero et al. \cite{MA} derived the existence of
a blow-up solution with $\int_{\mathbb{R}^2}u_0>\frac{8\pi}{\chi}$.
Contrastingly, any solution with $\int_{\mathbb{R}^2}u_0<\frac{8\pi}{\chi}$ exists
globally in time under additional conditions of $u_0$ \cite{VC}.
After this, Mizoguchi \cite{NM} proved that the system has global-in-time
solutions under no extra conditions.
When $N>2$, the existence of global solutions was deduced by Corrias and
Perthame \cite{LC} and the decay rates and asymptotic profiles of bounded
solutions have been given in \cite{TN3,TN1}.

The above mentioned references indicate that the solution of the classical
Keller-Segel system may occur the blow-up phenomenon.
The addition of a source term in the first equation of the model \eqref{1.2}
can prevent blow-up of solutions.
The common form of the source term is logistic growth which realistically
describes the dynamic response of the population.
The Keller-Segel system with Neumann initial-boundary value condition reads as
\begin{align} \label{1.3}\allowdisplaybreaks
 \left\{
    \begin{array}{llll}
        \displaystyle u_t=\Delta u-\chi\nabla\cdot(u\nabla v)+au-bu^\gamma, &&x \in\Omega, t>0,
        \\
        \displaystyle v_t=\Delta v+u-v, &&x \in \Omega, t>0,
        \\
        \displaystyle \frac{\partial u}{\partial n}=\frac{\partial v}{\partial n}=0,
        &&x\in\partial\Omega.
    \end{array}
 \right.
\end{align}
When $\gamma=2$, Osaki et al. \cite{KO2} derived the existence of the global solution
by establishing a prior estimate of the local solution for $N=1$.
When the space dimension $N=2$, the occurrence of blow-up is completely prevented
as long as a logistic source present and no critical mass blow-up is allowed
in this case \cite{KO3,TX1}.
In higher dimensions, the properly strong logistic damping can inhibit
the blow-up of solutions for \eqref{1.3}.
Winkler \cite{MW3} proved that the system \eqref{1.3} possesses
a unique global-in-time classical solution with $b>b_0$ sufficiently large.
Xiang \cite{TX2} extended the result of \cite{MW3} to obtain the
explicit form of logistic damping rate $b_0$. In other words, the problem of
how strong logistic damping is needed to prevent blow-up effectively resolved
in \cite{TX2}.
Moreover, the global asymptotic stability of the nontrivial equilibria
$(\frac{a}{b},\frac{a}{b})$ is given in \cite{XH,MW4}.
When $\gamma>1$, the system \eqref{1.3} possesses at least one global weak
solution under the condition that $\gamma>2-\frac{1}{N}$ for $N\geq2$ in \cite{GV1,GV2}
and the range of $\gamma$ is extended to $\gamma>\frac{2N+4}{N+4}$ 
for any space dimension in \cite{MW5}.

When the domain $\Omega$ is the entire space $\mathbb{R}^N$, consider
the following Cauchy problem for chemotaxis system with the logistic-type
source
\begin{align} \label{1.1}\allowdisplaybreaks
 \left\{
    \begin{array}{llll}
        \displaystyle u_t=\Delta u-\chi \nabla\cdot(u\nabla v )
        +au-bu^{\gamma}, &&x \in \mathbb{R}^N,t>0,
        \\
        \displaystyle v_t=\Delta v+\mu u-\lambda v, &&x \in \mathbb{R}^N,t>0,\\
        \displaystyle u(x,0)=u_0(x)\geq0, v(x,0)=v_0(x)\geq0, &&x \in \mathbb{R}^N,
    \end{array}
 \right.
\end{align}
where $N\geq2$, $\lambda>0$ describes the degradation rate of the chemical substance.
The term $\mu u$ represents the mobile species that produces the chemical
substance with $\mu>0$.
When $\gamma=2$, with simple calculations, it yields that
\begin{align*}
    \frac{d}{dt}\left(\frac{1}{2}|\nabla v|^2+\frac{1}{\chi}u\right)
    \leq\Delta\left(\frac{1}{2}|\nabla v|^2+\frac{1}{\chi}u\right)
    -|\nabla v|^2-\left(\frac{b}{\chi}-\frac{N}{4}\right)u^2+\frac{a}{\chi}u.
\end{align*}
The function $z=\frac{1}{2}|\nabla v|^2+\frac{1}{\chi}u$ satisfies
the following scalar parabolic inequality
\begin{align*}
    z_t\leq\Delta z-z+C
\end{align*}
with some positive $C$. For a bounded domain, the maximum principle
and the convexity of $\Omega$ ensure the boundedness of $z$,
which implies that both $u$ and $v$ are also bounded \cite{MW3}.
However, the approaches discussed in \cite{MW3} and other
references for bounded domain are not directly applicable to
the entire space. In recent years, there have been several papers on
the Cauchy problem \eqref{1.1} for $\gamma=2$.
Shen and Xue \cite{WS} proved that the global existence and
persistence phenomena of classical solution for every
nonnegative bounded and uniformly continuous initial data
if $b>\frac{N\mu\chi}{4}$ and under the assumption
$\lambda>\frac{a}{2}$, the constant equilibria
$(\frac{a}{b},\frac{\mu a}{\lambda b})$ is asymptotically
stable in the sense that
\begin{align*}
    \lim_{t\rightarrow\infty}[\|u(x,t)-\frac{a}{b}\|
    _{L^\infty(\mathbb{R}^N)}+\|v(x,t)-\frac{\mu a}{\lambda b}\|
    _{L^\infty(\mathbb{R}^N)}]=0.
\end{align*}
Similarly, if $b$ is large enough, the global-in-time boundedness
of nonnegative classical solution has been obtained in \cite{YN}
by developing localized estimations. In contrast to \cite{WS},
the upper bounds of solutions are given explicitly for small
$\chi$ \cite{DJXM}. Moreover, the assumption $\lambda>\frac{a}{2}$
in \cite{DJXM} is removed. We can conclude that the global boundedness
and asymptotic behavior of the global classical solution for the system
\eqref{1.1} hold when $\gamma=2$. However, it has not been determined
the well-posedness of the chemotaxis model \eqref{1.1} with
$1<\gamma<2$. Inspired by \cite{DJXM}, we consider the parabolic-parabolic
Keller-Segel system \eqref{1.1} with $\gamma>1$.
In the present paper, we extend the result in \cite{DJXM} to the case $\gamma>1$.
Replacing the analytic semigroup $\{e^{t(\Delta-a I)}\}_{t>0}$ by
$\{e^{t(\Delta-(\gamma-1)a I)}\}_{t>0}$, we obtain the global boundedness
of solutions and the explicit upper bounds of solutions
are as follows
\begin{align*}
    \|u(t)\|_{L^\infty(\mathbb{R}^N)}\leq e^{-(\gamma-1)at}\|u_0\|_{L^\infty(\mathbb{R^N})}
    +\frac{3}{2}\left(\frac{a}{b}\right)^{\frac{1}{\gamma-1}}.
\end{align*}
It can establish the groundwork for the asymptotic stability
of the positive constant equilibria.
Furthermore, the presence of the source term causes difficulties
in investigating the asymptotic stability of solutions when $1<\gamma<2$,
which can be solved by finding a quadratic term that controls the source term.

The present paper is organized as follows.
In Section 2, we provide some preliminary materials required for the proof of
our main results.
In Section 3, we list all the results of this paper.
In Section 4, we study the global boundedness of the classical solutions of \eqref{1.1}
with given initial functions when $\gamma>1$.
In Section 5, we discuss the asymptotic behavior of global classical solutions
when $1<\gamma<2$.

\section{Preliminaries}
In this section, we introduce some notations and well-known results which will be used
in the other sections.

Let
\begin{equation*}
    C_{unif}^b(\mathbb{R}^N)=\{u\in C(\mathbb{R}^N)|u(x)\text{ is uniformly continuous in }
    x\in\mathbb{R}^N \text{and} \sup_{x\in\mathbb{R}^N}|u(x)|<\infty\}
\end{equation*}
equipped with the norm $\|u\|_{L^\infty(\mathbb{R}^N)}=\sup_{x\in\mathbb{R}^N}|u(x)|$, and
\begin{equation*}
    C_{unif}^{b,1}(\mathbb{R}^N)=\{u\in C_{unif}^b(\mathbb{R}^N)|\partial_{x_i}u\in
    C_{unif}^b(\mathbb{R}^N), i=1,2,\cdots,N\}
\end{equation*}
equipped with the norm $\|u\|_{C_{unif}^{b,1}(\mathbb{R}^N)}=\|u\|_{L^\infty(\mathbb{R}^N)}
+\Sigma_{i=1}^{N}\|\partial_{x_i}u\|_{L^\infty(\mathbb{R}^N)}$ and
\begin{equation*}
    C_{unif}^{b,2}(\mathbb{R}^N)=\{u\in C_{unif}^{b,1}(\mathbb{R}^N)|\partial_{x_ix_j}u\in
    C_{unif}^b(\mathbb{R}^N), i=1,2,\cdots,N\}.
\end{equation*}
For given $0<\nu<1$, let
\begin{equation*}
    C_{unif}^{b,\nu}(\mathbb{R}^N)=\{u\in C_{unif}^b(\mathbb{R}^N)|\sup_{x,y\in\mathbb{R}^N,x\neq y}
    \frac{|u(x)-u(y)|}{|x-y|^\nu}<\infty\}
\end{equation*}
with the norm $\|u\|_{C_{unif}^{b,\nu}(\mathbb{R}^N)}=\sup_{x\in\mathbb{R}^N}|u(x)|
+\sup_{x,y\in\mathbb{R}^N,x\neq y}\frac{|u(x)-u(y)|}{|x-y|^\nu}$.
For $0<\theta<1$, let
\begin{align*}
    &C^\theta((t_1,t_2),C_{unif}^{b,\nu}(\mathbb{R}^N))\\
    &=\{u(\cdot)\in C((t_1,t_2),C_{unif}^{b,\nu}(\mathbb{R}^N))|u(t)\text{ is locally H\"{o}lder
    continuous with exponent } \theta\}.
\end{align*}

For $\rho>0$, let $\{e^{t(\Delta-\rho I)}\}_{t>0}$ denote the analytic semigroup
generated by $\Delta-\rho I$ on $C_{unif}^{b}(\mathbb{R}^N)$, then we have the
following basic statements.

\begin{lemma}(\cite{DH0,AP})\label{le2.1}
For every $t>0$,
\begin{align*}
    &\|e^{t(\Delta-\rho I)}u\|_{L^\infty(\mathbb{R}^N)}\leq
    e^{-\rho t}\|u\|_{L^\infty(\mathbb{R}^N)},\\
    &\|\nabla e^{t(\Delta-\rho I)}u\|_{L^\infty(\mathbb{R}^N)}\leq
    C_N t^{-\frac{1}{2}} e^{-\rho t}\|u\|_{L^\infty(\mathbb{R}^N)},
\end{align*}
holds for all $u\in L^\infty(\mathbb{R}^N)$, where $C_N>0$ is a constant depending only on $N$.
\end{lemma}

\begin{lemma}(Lemma 3.2 in \cite{RBS})\label{le2.2}
For every $t>0$, the operator $e^{t(\Delta-\rho I)}\nabla\cdot$ has a unique bounded extension
on $(C_{unif}^b(\mathbb{R}^N))^N$ satisfying
\begin{align*}
    \|e^{t(\Delta-\rho I)}\nabla\cdot u\|_{L^\infty(\mathbb{R}^N))}\leq\frac{N}{\sqrt{\pi}}
    t^{-\frac{1}{2}}e^{-\rho t}\|u\|_{L^\infty(\mathbb{R}^N)}, \quad \forall u\in (C_{unif}^b(\mathbb{R}^N))^N.
\end{align*}
\end{lemma}

\section{Statements of the main results}
The main aim of this paper is to investigate the global boundedness and asymptotic behavior
of a unique classical solution $(u(x,t),v(x,t))$ with initial function
$(u_0,v_0)\in C_{unif}^{b}(\mathbb{R}^N)\times C_{unif}^{b,1}(\mathbb{R}^N)$.
Motivated mainly from in \cite{WS}, we first obtain the local existence of classical
solutions of the system \eqref{1.1} by the contraction mapping theorem when $\gamma>1$.
To avoid repetition, we omit giving details on this here.

\begin{proposition} \label{pro4.1}
Assume $0\leq u_0\in C_{unif}^b(\mathbb{R}^N)$, $0\leq v_0\in C_{unif}^{b,1}(\mathbb{R}^N)$
and $\gamma\in(1,\infty)$.
There exists $T_{max}\in (0,\infty]$ such that the system \eqref{1.1} has a unique nonnegative
classical solution $(u,v)$ satisfying
\begin{align*}
    u\in C([0,T_{max});C_{unif}^b(\mathbb{R}^N))\cap C^1((0,T_{max});C_{unif}^b(\mathbb{R}^N)),\\
    v\in C([0,T_{max});C_{unif}^{b,1}(\mathbb{R}^N))\cap C^1((0,T_{max});C_{unif}^{b,1}(\mathbb{R}^N)),\\
    u,v,\partial_{x_i}u,\partial_{x_i}v,\partial_{x_i x_j}^2 u,\partial_{x_i x_j}^2 v,\partial_{t}u,
    \partial_{t}v\in C^\theta([t_0,T_{max});C_{unif}^{b,\nu}(\mathbb{R}^N))
\end{align*}
for all $i,j=1,2,\cdots,N$, $0<t_0<T_{max}$, $0<\theta\ll 1$, $0<\nu\ll 1$.
Furthermore, if $T_{max}<\infty$, then
\begin{align}\label{3.1}
    \lim_{t\rightarrow T_{max}^-}(\|u(\cdot,t)\|_{L^\infty(\mathbb{R}^N)}+\|v(\cdot,t)\|_{W^{1,\infty}
    (\mathbb{R}^N)})=\infty.
\end{align}
\end{proposition}

The first theorem concerns the global existence of the classical solution with initial function
$(u_0,v_0)\in C_{unif}^{b}(\mathbb{R}^N)\times C_{unif}^{b,1}(\mathbb{R}^N)$ when $\gamma>1$.

\begin{theorem}\label{th3.2}
Assume $0\leq u_0\in C_{unif}^b(\mathbb{R}^N)$, $0\leq v_0\in C_{unif}^{b,1}(\mathbb{R}^N)$
and $\gamma\in(1,\infty)$.
There exists $\chi_0>0$ depending on $a,b,\mu,\lambda,\gamma,N,\|u_0\|_{L^\infty(\mathbb{R}^N)},
\|\nabla v_0\|_{L^\infty(\mathbb{R}^N)}$ such that for all $\chi<\chi_0$,
the system \eqref{1.1} has a unique nonnegative classical global-in-time solution
\begin{align*}
    u\in C^0([0,\infty);L^\infty(\mathbb{R}^N))\cap C^{2,1}(\mathbb{R}^N,(0,\infty)),\\
    v\in C^0([0,\infty);W^{1,\infty}(\mathbb{R}^N))\cap C^{2,1}(\mathbb{R}^N,(0,\infty)),
\end{align*}
which enjoys global-in-time boundedness, that is, there exits a positive constant $C>0$ such that
for all $t>0$,
\begin{align*}
    \|u(\cdot,t)\|_{L^\infty(\mathbb{R}^N)}+\|v(\cdot,t)\|_{W^{1,\infty}(\mathbb{R}^N)}\leq C.
\end{align*}
\end{theorem}

The next proposition refers to the persistence of the global classical solution with strictly positive
initial $u_0$ when $\gamma>1$.

\begin{proposition}\label{pro5.1}
Suppose that $\gamma>1$ and $\chi<\chi_0$, then there exist $m>0$ and $M>0$
such that for any $u_0\in C_{unif}^b(\mathbb{R}^N)$,
$0\leq v_0\in C_{unif}^{b,1}(\mathbb{R}^N)$ with
$\inf_{x\in\mathbb{R}^N}u_0>0$, there exists $T$ such that
\begin{align*}
    m\leq u(x,t)\leq M
\end{align*}
for any $x\in\mathbb{R}^N$, $t\geq T$.
\end{proposition}

\begin{remark}
In Theorem 1.3 of \cite{WS}, the persistence of classical solutions of \eqref{1.1}
is obtained under the assumption $b>\frac{N\mu\chi}{4}$ for the case $\gamma=2$.
In this paper, we get Proposition \ref{pro5.1} by using the boundedness of
global classical solutions in Theorem \ref{th3.2} under the condition $\chi<\chi_0$.
Similar to the proof of Lemma 4.1 in \cite{WS}, using the global boundedness of
the classical solutions, we find that there are $M>0$ and $T_0>1$ such that for any $t\geq T_0$,
$\|u\|_{L^\infty(\mathbb{R}^N)}\leq M$, $\|v\|_{L^\infty(\mathbb{R}^N)}\leq M$,
$\|\nabla v\|_{L^\infty(\mathbb{R}^N)}\leq M$ and $\|\Delta v\|_{L^\infty(\mathbb{R}^N)}\leq M$
hold. Then we obtain Proposition \ref{pro5.1} by using the above estimates.
Different from Step 3 of the proof of Theorem 1.3 in \cite{WS},
we may assume that $\frac{3a}{4}-b\varepsilon_0^{\gamma-1}\geq\frac{a}{2}$.
From Step 1 of the proof in \cite{WS}, there is $\varepsilon_0>0$ such that for
$t_{1}+T(\varepsilon)+1\leq t<t_{2}\leq\infty$,
\begin{align*}
    u_{t}&=\Delta u-\chi\nabla v\cdot\nabla u+u(a-\chi\Delta v-bu^{\gamma-1})\\
    &\geq\Delta u+q(x,t)\cdot\nabla u+\frac{a}{2}u
\end{align*}
where $q(x,t)=-\chi\nabla v(x,t)$.
Since the rest of the proof is similar to that of \cite{WS}, we omit it here.
\end{remark}

The last theorem focuses on the asymptotic behavior of the global classical solution with $1<\gamma<2$.

\begin{theorem}\label{th3.4}
Assume that $u_0\in C_{\text {unif }}^{b}(\mathbb{R}^{N})$,
$0\leq v_{0}\in C_{\text{unif}}^{b, 1}(\mathbb{R}^{N})$
, $\inf _{x\in\mathbb{R}^{N}}u_{0}>0$ and $\gamma\in(1,2)$.
There exist $T_{m}>0$, $C>0$, $\sigma\in(0,\min\{\frac{a}{2},\frac{\lambda}{\gamma-1}\})$,
and $\chi_{*}>0$ depending on $\inf_{x\in\mathbb{R}^{N}} u_{0},a,b,\mu,\gamma$,
$\lambda,N,\left\|u_{0}\right\|_{L^{\infty}\left(\mathbb{R}^{N}\right)},
\left\|\nabla v_{0}\right\|_{L^{\infty}\left(\mathbb{R}^{N}\right)}$, such that for $\chi<\chi_*$,
the nonnegative classical solution of problem \eqref{1.1} satisfies
\begin{align*}
    \left\|u(t)-\left(\frac{a}{b}\right)^{\frac{1}{\gamma-1}}\right\|_{L^{\infty}\left(\mathbb{R}^{N}\right)}
    \leq C e^{-\sigma(\gamma-1) t}, \quad
    \left\|v(t)-\frac{\mu}{\lambda}\left(\frac{a}{b}\right)^{\frac{1}{\gamma-1}}\right\|_{L^{\infty}
    \left(\mathbb{R}^{N}\right)} \leq C e^{-\sigma(\gamma-1) t}
\end{align*}
for all $t>T_{m}$.
\end{theorem}

\section{Global boundedness of classical solutions}\label{S3}
This section is dedicated to the investigation of the local and global boundedness of
classical solutions for the problem \eqref{1.1} with nonnegative initial functions.

We consider the following Cauchy problem
\begin{align} \label{4.2}\allowdisplaybreaks
 \left\{
    \begin{array}{llll}
        \displaystyle \frac{\partial w}{\partial t } =\Delta w-\rho w+f, &&x \in \mathbb{R}^N,t>0,
        \\
        \displaystyle w(x,0)=w_0(x), &&x \in \mathbb{R}^N.
    \end{array}
 \right.
\end{align}
In order to state the results, we give the definition of a mild solution to \eqref{4.2}.
\begin{definition}
Let $w_0\in C_{unif}^{b}(\mathbb{R}^N)$ and $f\in L^1((0,+\infty);C_{unif}^{b}(\mathbb{R}^N))$.
The function $w(x,t)\in C([0,+\infty);C_{unif}^{b}(\mathbb{R}^N))$ given by
\begin{align*}
    w(t)=e^{t(\Delta-\rho I)}w_0(x)+\int_{0}^{t}e^{(t-s)(\Delta-\rho I)}f(s)ds
\end{align*}
for $0\leq t<+\infty$ is the mild solution of \eqref{4.2} on $[0,+\infty)$, we have
\begin{align*}
    e^{t(\Delta-\rho I)}w(x)=\int_{\mathbb{R}^N}\frac{1}{(4\pi t)^{\frac{N}{2}}}
    e^{-\frac{|x-y|^2}{4t}-\rho t}w(y)dy.
\end{align*}
\end{definition}

{\bf Proof of Theorem \ref{th3.2}.}
To derive the global boundedness of classical solutions, we rewrite \eqref{1.1} as
\begin{align*} \allowdisplaybreaks
 \left\{
    \begin{array}{llll}
        \displaystyle \frac{\partial u}{\partial t } =(\Delta-(\gamma-1)a) u-\chi\nabla\cdot(u\nabla v )
        +\gamma au-bu^\gamma,     &&x \in \mathbb{R}^N,t>0,
        \\
        \displaystyle \frac{\partial v}{\partial t } =\Delta v+\mu u-\lambda v, &&x \in \mathbb{R}^N,t>0
    \end{array}
 \right.
\end{align*}
with $u(x,0)=u_{0}(x)\in C_{unif}^{b}(\mathbb{R}^N)$, and $v(x,0)=v_{0}(x)\in C_{unif}^{b,1}(\mathbb{R}^N)$
for $x \in \mathbb{R}^N$. It is easy to verify that for $\gamma>1$, we have
\begin{align}\label{3.3}
    u(t)&=e^{t(\Delta-(\gamma-1)aI)}u_0-\chi\int_{0}^{t}e^{s(\Delta-(\gamma-1)aI)}\nabla\cdot(u\nabla v)
    (t-s)ds\nonumber\\
    &~~~+\int_{0}^{t}e^{s(\Delta-(\gamma-1)aI)}(\gamma au(t-s)-bu^\gamma(t-s)) ds,\\
    v(t)&=e^{t(\Delta-\lambda I)}v_0+\mu\int_{0}^{t}e^{s(\Delta-\lambda I)}u(t-s)ds. \label{3.4}
\end{align}
Define
\begin{align*}
    \hat{T}=\sup\{T_* \in(0,T_{max});\|u(t)\|_{L^\infty(\mathbb{R}^N)}\leq \bar{C},  t\in[0,T_*)\},
\end{align*}
where
\begin{align}\label{3.6}
    \bar{C}=2\|u_0\|_{L^\infty(\mathbb{R}^N)}+3\left(\frac{a}{b}\right)^{\frac{1}{\gamma-1}}.
\end{align}
Notice that $\|u_0\|_{L^\infty(\mathbb{R}^N)}<\bar{C}$ and $\hat{T}>0$ is well-defined.

{\bf Claim 1.}
\textit{For any} $t\in(0,\hat{T})$, \textit{we have}
\begin{align*}
    \|v(t)\|_{L^\infty(\mathbb{R}^N)}\leq\|v_0\|_{L^\infty(\mathbb{R}^N)}+\frac{\mu}{\lambda}\bar{C}
\end{align*}
\textit{and}
\begin{align*}
    \|\nabla v(t)\|_{L^\infty(\mathbb{R}^N)}\leq\|\nabla v_0\|_{L^\infty(\mathbb{R}^N)}+\mu \bar{C}C_NC_1.
\end{align*}

It follows from \eqref{3.4} and Lemma \ref{le2.1} that
\begin{align*}
    \|v(t)\|_{L^\infty(\mathbb{R}^N)}
    &\leq\|e^{t(\Delta-\lambda I)}v_0\|_{L^\infty(\mathbb{R}^N)}
    +\mu\int_{0}^{t}\|e^{s(\Delta-\lambda I)}u(t-s)\|_{L^\infty(\mathbb{R}^N)}ds\\
    &\leq e^{-\lambda t}\|v_0\|_{L^\infty(\mathbb{R}^N)}+\mu\int_{0}^{t}e^{-\lambda s}\|u(t-s)\|_{L^\infty(\mathbb{R}^N)}ds\\
    &\leq\|v_0\|_{L^\infty(\mathbb{R}^N)}+\frac{\mu}{\lambda}\bar{C}
\end{align*}
for all $t\in(0,\hat{T})$, where $\bar{C}$ is defined in \eqref{3.6}.

It is easy to see that $\nabla v(t)=\nabla(e^{t(\Delta-\lambda I)}v_0)
+\mu\int_{0}^{t}\nabla(e^{s(\Delta-\lambda I)}u(t-s))ds$ for all $t\in(0,\hat{T})$.
According to $\nabla(e^{t(\Delta-\lambda I)}v_0)=e^{t(\Delta-\lambda I)}\nabla v_0$, we get
\begin{align*}
    \|\nabla v(t)\|_{L^\infty(\mathbb{R}^N)}
    &\leq\|e^{t(\Delta-\lambda I)}\nabla v_0\|_{L^\infty(\mathbb{R}^N)}
    +\mu\int_{0}^{t}\|\nabla(e^{s(\Delta-\lambda I)}u(t-s))\|_{L^\infty(\mathbb{R}^N)}ds\\
    &\leq e^{-\lambda t}\|\nabla v_0\|_{L^\infty(\mathbb{R}^N)}+\mu C_N\int_{0}^{t}s^{-\frac{1}{2}}
    e^{-\lambda s}\|u(t-s)\|_{L^\infty(\mathbb{R}^N)}ds\\
    &\leq \|\nabla v_0\|_{L^\infty(\mathbb{R}^N)}+\mu \bar{C}C_NC_1,
\end{align*}
where $C_N$ is the constant in Lemma \ref{le2.1} and $C_1=\int_{0}^{+\infty}s^{-\frac{1}{2}}e^{-\lambda s}ds<+\infty$.

{\bf Claim 2.}
\textit{For any} $t\in(0,\hat{T})$, \textit{we have}
$$
\|u(t)\|_{L^\infty(\mathbb{R}^N)}\leq\frac{\bar{C}}{2}.
$$

Since $u\geq0$, thus for $\gamma>1$, we have
\begin{align*}
    \gamma au-bu^\gamma\leq a(\gamma-1)
    \left(\frac{a}{b}\right)^{\frac{1}{\gamma-1}}.
\end{align*}
This combined with \eqref{3.3} yield that
\begin{align*}
    u(t)&=e^{t(\Delta-(\gamma-1)aI)}u_0-\chi\int_{0}^{t}e^{s(\Delta-(\gamma-1)aI)}\nabla\cdot(u\nabla v)(t-s)ds\\
    &~~~+\int_{0}^{t}e^{s(\Delta-(\gamma-1)aI)}(\gamma au(t-s)-bu^\gamma(t-s)) ds,\\
    &\leq e^{t(\Delta-(\gamma-1)aI)}u_0-\chi\int_{0}^{t}e^{s(\Delta-(\gamma-1)aI)}\nabla\cdot(u\nabla v)(t-s)ds\\
    &~~~+\int_{0}^{t}e^{s(\Delta-(\gamma-1)aI)}a(\gamma-1)\left(\frac{a}{b}\right)^{\frac{1}{\gamma-1}}ds.
\end{align*}
By Lemma \ref{le2.1} and Lemma \ref{le2.2}, it can be derived that
\begin{align*}
    \|u(t)\|_{L^\infty(\mathbb{R}^N)}&\leq \|e^{t(\Delta-(\gamma-1)aI)}u_0\|_{L^\infty(\mathbb{R^N})}
    +\chi\int_{0}^{t}\|e^{s(\Delta-(\gamma-1)aI)}\nabla\cdot(u\nabla v)(t-s)\|_{L^\infty(\mathbb{R^N})}ds\\
    &~~~+\int_{0}^{t}\|e^{s(\Delta-(\gamma-1)aI)}a(\gamma-1)\left(\frac{a}{b}\right)^{\frac{1}{\gamma-1}}
    \|_{L^\infty(\mathbb{R^N})}ds\\
    &\leq e^{-(\gamma-1)at}\|u_0\|_{L^\infty(\mathbb{R^N})}+\frac{\chi N}{\sqrt{\pi}}
    \sup_{t\geq0}\|u(t)\|_{L^\infty(\mathbb{R}^N)}\sup_{t\geq0}\|\nabla v(t)\|_{L^\infty(\mathbb{R}^N)}
    \int_{0}^{t}s^{-\frac{1}{2}}e^{-(\gamma-1)as}ds\\
    &~~~+a(\gamma-1)\left(\frac{a}{b}\right)^{\frac{1}{\gamma-1}}\int_{0}^{t}e^{-(\gamma-1)as}ds\\
    &\leq e^{-(\gamma-1)at}\|u_0\|_{L^\infty(\mathbb{R^N})}+\chi \bar{C}C_2(\|\nabla v_0\|_{L^\infty(\mathbb{R}^N)}
    +\mu \bar{C}C_NC_1)+\left(\frac{a}{b}\right)^{\frac{1}{\gamma-1}}
\end{align*}
for all $t\in(0,\hat{T})$ with $C_2=\frac{N}{\sqrt{\pi}}\int_{0}^{+\infty}s^{-\frac{1}{2}}e^{-(\gamma-1)as}ds<+\infty$.
Using the definition of $\bar{C}$ in \eqref{3.6}, for
\begin{align}\label{3.7}
    \chi<\chi_0:=\frac{\left(\frac{a}{b}\right)^{\frac{1}{\gamma-1}}}
    {2\bar{C}C_2(\|\nabla v_0\|_{L^\infty(\mathbb{R}^N)}+\mu \bar{C}C_NC_1)},
\end{align}
we can get
\begin{align}\label{3.8}
    \|u(t)\|_{L^\infty(\mathbb{R}^N)}\leq e^{-(\gamma-1)at}\|u_0\|_{L^\infty(\mathbb{R^N})}
    +\frac{3}{2}\left(\frac{a}{b}\right)^{\frac{1}{\gamma-1}}
    \leq\frac{\bar{C}}{2}.
\end{align}
Hence $\|u(t)\|_{L^\infty(\mathbb{R}^N)}\leq\frac{\bar{C}}{2}$ holds for all $t\in(0,\hat{T})$.

Now, we show that $\hat{T}=T_{max}$ by contraction.
Indeed, suppose that  $\hat{T}<T_{max}$.
Then it follows from Claim 2 that
$\|u(\cdot,t)\|_{L^{\infty}(\mathbb{R}^{N})}\leq\frac{\bar{C}}{2}$ for all $t \in(0, \hat{T})$,
which contradicts the definition of $\hat{T}$.

Next, we prove that $T_{max}=+\infty$.
Indeed, if $T_{max}<+\infty$, we conclude from $\hat{T}=T_{max}$, Claim 1
and Claim 2 that
\begin{align*}
    \lim _{t \rightarrow T_{max}^{-}}\left(\|u(\cdot, t)\|_{L^{\infty}\left(\mathbb{R}^{N}\right)}
    +\|v(\cdot, t)\|_{W^{1, \infty}\left(\mathbb{R}^{N}\right)}\right)<\infty,
\end{align*}
which contradicts with \eqref{3.1} in Proposition \ref{pro4.1}.
Thus, we have $T_{max}=+\infty$, which finishes the proof of Theorem \ref{th3.2}.
$\hfill\Box$

\section{Asymptotic stability of constant equilibria}
The goal of this section is to prove the asymptotic behavior of global bounded
classical solution of \eqref{1.1}.
The innovative ideas what we shall employ are mainly from \cite{DJXM}.

Define
\begin{align*}
    \tilde{u}=u-\left(\frac{a}{b}\right)^{\frac{1}{\gamma-1}},
    \quad \tilde{v}=v-\frac{\mu}{\lambda}\left(\frac{a}{b}\right)^{\frac{1}{\gamma-1}},
\end{align*}
and focus on the following system
\begin{align*}\allowdisplaybreaks
 \left\{
    \begin{array}{llll}
        \displaystyle \tilde{u}_t&=(\Delta-(\gamma-1)aI)\tilde{u}-\chi \nabla\cdot(\tilde{u}\nabla \tilde{v})
        -\chi\left(\frac{a}{b}\right)^{\frac{1}{\gamma-1}}\Delta\tilde{v}\\
        &~~~+a(\gamma \tilde{u}+\left(\frac{a}{b}\right)^{\frac{1}{\gamma-1}})
        -b(\tilde{u}+\left(\frac{a}{b}\right)^{\frac{1}{\gamma-1}})^\gamma,
        &&x \in \mathbb{R}^N,t>0,
        \\
        \displaystyle \tilde{v}_t&=\Delta \tilde{v}+\mu\tilde{u}-\lambda\tilde{v},
        &&x \in \mathbb{R}^N,t>0
    \end{array}
 \right.
\end{align*}
with $\tilde{u}(x,0)=u_0(x)-\left(\frac{a}{b}\right)^{\frac{1}{\gamma-1}}$
and $\tilde{v}(x,0)=v_0(x)-\frac{\mu}{\lambda}\left(\frac{a}{b}\right)^{\frac{1}{\gamma-1}}$
for $x\in\mathbb{R}^N$.
Then for any $T>0$, we have $\tilde{u}$ and $\tilde{v}$ satisfy
\begin{align}
    \tilde{u}(t)&=e^{(t-T)(\Delta-(\gamma-1)aI)}\tilde{u}(T)
    -\chi\left(\frac{a}{b}\right)^{\frac{1}{\gamma-1}}\int_T^t e^{(s-T)(\Delta-(\gamma-1)aI)}
    \Delta\tilde{v}(t-s+T)ds\nonumber\\
    &~~~-\chi\int_T^t e^{(s-T)(\Delta-(\gamma-1)aI)}\nabla\cdot
    (\tilde{u}(t-s+T)\nabla \tilde{v}(t-s+T))ds\nonumber\\
    &~~~+\int_T^t e^{(s-T)(\Delta-(\gamma-1)aI)}\left[a(\gamma\tilde{u}(t-s+T)+
    \left(\frac{a}{b}\right)^{\frac{1}{\gamma-1}})
    -b(\tilde{u}(t-s+T)+\left(\frac{a}{b}\right)
    ^{\frac{1}{\gamma-1}})^\gamma\right] ds,\label{6.2}\\
    &\tilde{v}(t)=e^{(t-T)(\Delta-\lambda I)}\tilde{v}(T)+\mu\int_T^t e^{(s-T)(\Delta-\lambda I)}
    \tilde{u}(t-s+T)ds. \label{6.3}
\end{align}

\begin{lemma}\label{le4.4}
Under the assumption of Theorem \ref{th3.4}, there exists $T_{m}>0$ such that
\begin{align*}
    \|\tilde{u}(\cdot,t)\|_{L^{\infty}(\mathbb{R}^N)}
    <\frac{1}{4}\left(\frac{a}{b}\right)^{\frac{1}{\gamma-1}}
\end{align*}
for all $t\geq T_{m}$.
\end{lemma}

{\bf Proof.}
By Proposition \ref{pro5.1}, there exists
$\xi\in(0,\frac{1}{4}\left(\frac{a}{b}\right)^{\frac{1}{\gamma-1}})$ and
$T_0>0$ such that
\begin{align}\label{6.4}
    u(x,t)\geq\xi
\end{align}
for all $x\in \mathbb{R}^N$ and $t\in[T_0,+\infty)$.
In view of \eqref{3.8} and $\hat{T}=+\infty$, there exists $T_1>T_0$ such that
\begin{align*}
    \limsup_{t\rightarrow\infty}\|u(t)\|_{L^\infty(\mathbb{R}^N)}
    \leq\frac{3}{2}\left(\frac{a}{b}\right)^{\frac{1}{\gamma-1}}
    <2\left(\frac{a}{b}\right)^{\frac{1}{\gamma-1}}-\xi.
\end{align*}
holds for all $t\geq T_1$, which together with \eqref{6.4} imply
\begin{align*}
    \|\tilde{u}(t)\|_{L^\infty(\mathbb{R}^N)}
    \leq\left(\frac{a}{b}\right)^{\frac{1}{\gamma-1}}-\xi
\end{align*}
for all $t\geq T_1$.

Noticing
\begin{align*}
    \lim_{z\rightarrow 0}\frac{|a(\gamma z+
    \left(\frac{a}{b}\right)^{\frac{1}{\gamma-1}})-b(z
    +\left(\frac{a}{b}\right)^{\frac{1}{\gamma-1}})^\gamma|}
    {z^2}=\frac{a}{2}\gamma(\gamma-1)
    \left(\frac{b}{a}\right)^{\frac{1}{\gamma-1}},
\end{align*}
there exists $A>0$ such that for any $z$ satisfying
$|z|<A$, we have
\begin{align*}
    |a(\gamma z+
    \left(\frac{a}{b}\right)^{\frac{1}{\gamma-1}})-b(z
    +\left(\frac{a}{b}\right)^{\frac{1}{\gamma-1}})^\gamma|
    \leq a(\gamma-1)\left(\frac{b}{a}\right)^{\frac{1}{\gamma-1}}z^2.
\end{align*}
According to the strong maximum principle, $u$ is positive in
$\mathbb{R}^N\times(0,\infty)$. Then we have
$\tilde{u}>-\left(\frac{b}{a}\right)^{\frac{1}{\gamma-1}}$.
It can easily be verified that
\begin{align}\label{6.4.5}
    |a(\gamma\tilde{u}(t-s+T_1)+
    \left(\frac{a}{b}\right)^{\frac{1}{\gamma-1}})-b(\tilde{u}(t-s+T_1)
    +\left(\frac{a}{b}\right)^{\frac{1}{\gamma-1}})^\gamma|
    \leq a(\gamma-1)\left(\frac{b}{a}\right)^{\frac{1}{\gamma-1}}\tilde{u}^2
\end{align}
for any $\tilde{u}>-\left(\frac{b}{a}\right)^{\frac{1}{\gamma-1}}$.
According to \eqref{6.2}, \eqref{6.4.5}, Lemma \ref{le2.1} and Lemma \ref{le2.2},
we have
\begin{align*}
    \|\tilde{u}(t)\|_{L^\infty(\mathbb{R}^N)}
    &\leq\|e^{(t-T_1)(\Delta-(\gamma-1)aI)}\tilde{u}(T_1)\|_{L^\infty(\mathbb{R}^N)}\\
    &~~~+\chi\left(\frac{a}{b}\right)^{\frac{1}{\gamma-1}}
    \int_{T_1}^t\|e^{(s-T_1)(\Delta-(\gamma-1)aI)}\Delta\tilde{v}(t-s+T_1)\|_{L^\infty(\mathbb{R}^N)}ds\\
    &~~~+\chi\int_{T_1}^t\|e^{(s-T_1)(\Delta-(\gamma-1)aI)}\nabla\cdot(\tilde{u}(t-s+T_1)
    \nabla\tilde{v}(t-s+T_1))\|_{L^\infty(\mathbb{R}^N)}ds\\
    &~~~+\int_{T_1}^t \|e^{(s-T_1)(\Delta-(\gamma-1)aI)}[a(\gamma\tilde{u}(t-s+T_1)+
    \left(\frac{a}{b}\right)^{\frac{1}{\gamma-1}})
    -b(\tilde{u}(t-s+T_1)+\left(\frac{a}{b}\right)
    ^{\frac{1}{\gamma-1}})^\gamma]\|_{L^\infty(\mathbb{R}^N)}ds\nonumber\\
    &\leq e^{-a(\gamma-1)(t-T_1)}\|\tilde{u}(T_1)\|_{L^\infty(\mathbb{R}^N)}\\
    &~~~+a(\gamma-1)\left(\frac{b}{a}\right)^{\frac{1}{\gamma-1}}
    \int_{T_1}^te^{-a(\gamma-1)(s-T_1)}\|\tilde{u}^2(t-s+T_1)\|_{L^\infty(\mathbb{R}^N)}ds\\
    &~~~+\frac{\chi N}{\sqrt{\pi}}\left(\frac{a}{b}\right)^{\frac{1}{\gamma-1}}
    \int_{T_1}^te^{-a(\gamma-1)(s-T_1)}(s-T_1)^{-\frac{1}{2}}\|\nabla\tilde{v}(t-s+T_1)
    \|_{L^\infty(\mathbb{R}^N)}ds\\
    &~~~+\frac{\chi N}{\sqrt{\pi}}\int_{T_1}^te^{-a(\gamma-1)(s-T_1)}(s-T_1)^{-\frac{1}{2}}
    \|\tilde{u}(t-s+T_1)\|_{L^\infty(\mathbb{R}^N)}
    \|\nabla\tilde{v}(t-s+T_1)\|_{L^\infty(\mathbb{R}^N)}ds\\
    &\leq e^{-a(\gamma-1)(t-T_1)}(\left(\frac{a}{b}\right)^{\frac{1}{\gamma-1}}-\xi)
    +\left(\frac{b}{a}\right)^{\frac{1}{\gamma-1}}
    (\left(\frac{a}{b}\right)^{\frac{1}{\gamma-1}}-\xi)^2\\
    &~~~+\frac{\chi N}{\sqrt{\pi}}(2\left(\frac{a}{b}\right)^{\frac{1}{\gamma-1}}-\xi)
    (\|\nabla v_0\|_{L^\infty(\mathbb{R}^N)}+\mu \bar{C}C_NC_1)
    \int_{T_1}^te^{-a(\gamma-1)(s-T_1)}(s-T_1)^{-\frac{1}{2}}ds,
\end{align*}
where $\int_{T_1}^te^{-a(\gamma-1)(s-T_1)}(s-T_1)^{-\frac{1}{2}}ds
\leq\int_{0}^{+\infty}e^{-a(\gamma-1)s}s^{-\frac{1}{2}}ds<+\infty$.
Let
\begin{align}\label{6.5}
    \chi_1=\frac{\sqrt{\pi}\xi}{4N(2\left(\frac{a}{b}\right)^{\frac{1}{\gamma-1}}-\xi)
    (\|\nabla v_0\|_{L^\infty(\mathbb{R}^N)}+\mu \bar{C}C_NC_1)
    \int_{0}^{+\infty}e^{-a(\gamma-1)s}s^{-\frac{1}{2}}ds},
\end{align}
there exists $T_2>0$ such that
\begin{align}\label{6.5.5}
    \|\tilde{u}(t)\|_{L^\infty(\mathbb{R}^N)}\leq\frac{\xi}{2}
    +\left(\frac{b}{a}\right)^{\frac{1}{\gamma-1}}(\left(\frac{a}{b}\right)
    ^{\frac{1}{\gamma-1}}-\xi)^2
\end{align}
for all $t\geq T_2$ with $\chi<\min\{\chi_0,\chi_1\}$, where $\chi_0$ defined in \eqref{3.7}.

Then there exists $T_2>0$ such that
\begin{align*}
    \frac{\xi}{2}
    +\left(\frac{b}{a}\right)^{\frac{1}{\gamma-1}}(\left(\frac{a}{b}\right)
    ^{\frac{1}{\gamma-1}}-l_0\xi)^2
    =\left(\frac{a}{b}\right)^{\frac{1}{\gamma-1}}-l_1\xi
\end{align*}
with $l_0=1, l_1=\frac{3}{2}-\left(\frac{b}{a}\right)^{\frac{1}{\gamma-1}}\xi>1$,
which upon this observation yields that there exists a sequence $\{l_n\}_{n=1}^{\infty}$ satisfying
\begin{align*}
    l_n=-\left(\frac{b}{a}\right)^{\frac{1}{\gamma-1}}l_{n-1}^2\xi
    +2l_{n-1}-\frac{1}{2}
\end{align*}
and a nondecreasing sequence $\{T_n\}_{n=1}^{\infty}$ such that
\begin{align*}
    \|\tilde{u}(t)\|_{L^\infty(\mathbb{R}^N)}
    \leq\frac{\xi}{2}+\left(\frac{b}{a}\right)^{\frac{1}{\gamma-1}}
    (\left(\frac{a}{b}\right)^{\frac{1}{\gamma-1}}-l_n\xi)^2
    =\left(\frac{a}{b}\right)^{\frac{1}{\gamma-1}}-l_{n+1}\xi
\end{align*}
for all $t\geq T_{n+2}$.
From $\left(\frac{a}{b}\right)^{\frac{1}{\gamma-1}}-l_{n}\xi-\frac{\xi}{2}\geq0$,
we have $l_n<\left(\frac{a}{b}\right)^{\frac{1}{\gamma-1}}\frac{1}{\xi}-\frac{1}{2}$
for all $n\geq1$.
Furthermore, we easily check that $l_n\geq l_{n-1}\geq\cdots\geq l_1>l_0=1$
 for all $n\geq1$.
Thus, the sequence $\{l_n\}_{n=1}^{\infty}$ converges such that
\begin{align*}
    \lim_{n\rightarrow\infty}l_n=\tilde{l}
    =\left(\frac{a}{b}\right)^{\frac{1}{\gamma-1}}\frac{1}{2\xi}
    \big(1+\sqrt{1-2\left(\frac{b}{a}\right)^{\frac{1}{\gamma-1}}\xi}\big).
\end{align*}
Then we have
\begin{align*}
    \left(\frac{a}{b}\right)^{\frac{1}{\gamma-1}}-\tilde{l}\xi
    =\frac{1}{2}\left(\frac{a}{b}\right)^{\frac{1}{\gamma-1}}
    (1-\sqrt{1-2\left(\frac{b}{a}\right)^{\frac{1}{\gamma-1}}\xi})
    \leq\frac{1}{2}\left(\frac{a}{b}\right)^{\frac{1}{\gamma-1}}(1-\frac{\sqrt{2}}{2}).
\end{align*}
Thus we conclude that there exists $m\in\mathbb{N}$ such that for any $t\geq T_{m}$,
\begin{align*}
    \|\tilde{u}(t)\|_{L^\infty(\mathbb{R}^N)}
    \leq\left(\frac{a}{b}\right)^{\frac{1}{\gamma-1}}-l_{m-1}\xi
    <\frac{1}{4}\left(\frac{a}{b}\right)^{\frac{1}{\gamma-1}}.
\end{align*}
$\hfill\Box$

{\bf Proof of Theorem \ref{th3.4}.}
Based on Lemma \ref{le4.4}, we can know that there exists
$\varepsilon\in(0,1)$ and $\sigma\in(0,\min\{\frac{a}{2},
\frac{\lambda}{\gamma-1}\})$
such that
\begin{align}\label{6.7}
    \|\tilde{u}(T_{m})\|_{L^\infty(\mathbb{R}^N)}\leq\frac{\varepsilon^2
    a^{\frac{2-\gamma}{\gamma-1}}(a-2\sigma)}{4b^{\frac{1}{\gamma-1}}}.
\end{align}
Let
\begin{align}\label{6.8}
    \tilde{T}=\sup\{T_*\in(T_{m},+\infty);\|\tilde{u}(T_{m})\|_{L^\infty(\mathbb{R}^N)}
    \leq \tilde{C}e^{-\sigma(\gamma-1)(t-T_{m})}, t\in[T_{m},T_*)\}
\end{align}
with $\tilde{C}=\frac{\varepsilon a^{\frac{2-\gamma}{\gamma-1}}(a-\sigma)}{2b^{\frac{1}{\gamma-1}}}$.
It follows from \eqref{6.7} that $\tilde{C}>\|\tilde{u}(T_{m})\|_{L^\infty(\mathbb{R}^N)}$ and
$\tilde{T}>T_{m}$ is well-defined.
Thus, the asymptotic stability of the positive constant equilibria
$(\left(\frac{a}{b}\right)^{\frac{1}{\gamma-1}},
\frac{\mu}{\lambda}\left(\frac{a}{b}\right)^{\frac{1}{\gamma-1}})$
for all $t\in(T_{m}, \tilde{T})$
can be derived from the following claims and
$\tilde{T}=+\infty$ can be proved by contradiction.

{\bf Claim 1.} \textit{If} $(u_0,v_0)\in C_{\text {unif }}^{b}(\mathbb{R}^{N})\times C_{\text {unif }}
^{b,1}(\mathbb{R}^{N})$, $\inf _{x\in\mathbb{R}^{N}}u_{0}>0$ \textit{and} $v_0\geq0$, \textit{we have}
\begin{align*}
    \|\tilde{v}(t)\|_{L^\infty(\mathbb{R}^N)}\leq
    e^{-\sigma(\gamma-1)(t-T_{m})}\left(\|v_0\|_{L^\infty(\mathbb{R}^N)}
    +\frac{\mu}{\lambda}\bar{C}+\frac{\mu}{\lambda-\sigma(\gamma-1)}\tilde{C}\right)
\end{align*}
\textit{and}
\begin{align*}
    \|\nabla\tilde{v}(t)\|_{L^\infty(\mathbb{R}^N)}\leq
    C_{\tilde{v}}e^{-\sigma(\gamma-1)(t-T_{m})}
\end{align*}
\textit{with}
\begin{align*}
    C_{\tilde{v}}=\|\nabla v_0\|_{L^\infty(\mathbb{R}^N)}+\mu \bar{C}C_NC_1+\mu \tilde{C}C_NC_3
\end{align*}
\textit{for all} $t\in(T_{m},\tilde{T})$.

Combining \eqref{6.3} with Lemma \ref{le2.1}, we deduce
\begin{align*}
     \|\tilde{v}(t)\|_{L^\infty(\mathbb{R}^N)}
     &\leq \|e^{(t-T_{m})(\Delta-\lambda I)}\tilde{v}(T_{m})\|_{L^\infty(\mathbb{R}^N)}
     +\mu\int_{T_{m}}^t\|e^{(s-T_{m})(\Delta-\lambda I)}\tilde{u}(t-s+T_{m})\|_{L^\infty(\mathbb{R}^N)}ds\\
     &\leq e^{-\lambda(t-T_{m})}\|\tilde{v}(T_{m})\|_{L^\infty(\mathbb{R}^N)}
     +\mu\int_{T_{m}}^te^{-\lambda(s-T_{m})}\|\tilde{u}(t-s+T_{m})\|_{L^\infty(\mathbb{R}^N)}ds\\
     &\leq e^{-\lambda(t-T_{m})}(\|v_0\|_{L^\infty(\mathbb{R}^N)}+\frac{\mu}{\lambda}\bar{C})
     +\mu \tilde{C}\int_{T_{m}}^te^{-\lambda(s-T_{m})}e^{-\sigma(\gamma-1)(t-s+T_{m})}ds\\
     &\leq e^{-\lambda(t-T_{m})}(\|v_0\|_{L^\infty(\mathbb{R}^N)}+\frac{\mu}{\lambda}\bar{C})
     +\mu \tilde{C}e^{-\sigma(\gamma-1) t}\int_0^{t-T_{m}}e^{-(\lambda-\sigma(\gamma-1))s}ds\\
     &\leq e^{-\sigma(\gamma-1)(t-T_{m})}(\|v_0\|_{L^\infty(\mathbb{R}^N)}+\frac{\mu}{\lambda}\bar{C}
     +\frac{\mu}{\lambda-\sigma(\gamma-1)}\tilde{C}).
\end{align*}
It is easy to see that $\nabla\tilde{v}(t)=e^{(t-T_{m})(\Delta-\lambda I)}\nabla\tilde{v}(T_{m})
+\mu\int_{T_{m}}^{t}\nabla(e^{(s-T_{m})(\Delta-\lambda I)}\tilde{u}(t-s+T_{m}))ds$.
Applying Lemma \ref{le2.1} and $\lambda-\sigma(\gamma-1)>0$, we derive
\begin{align*}
    \|\nabla\tilde{v}(t)\|_{L^\infty(\mathbb{R}^N)}
    &\leq \|e^{(t-T_{m})(\Delta-\lambda I)}\nabla\tilde{v}(T_{m})\|_{L^\infty(\mathbb{R}^N)}
    +\mu\int_{T_{m}}^t\|\nabla(e^{(s-T_{m})(\Delta-\lambda I)}\tilde{u}(t-s+T_{m}))\|_{L^\infty(\mathbb{R}^N)}ds\\
    &\leq e^{-\lambda(t-T_{m})}\|\nabla\tilde{v}(T_{m})\|_{L^\infty(\mathbb{R}^N)}
    +\mu C_N\int_{T_{m}}^t e^{-\lambda(s-T_{m})}(s-T_{m})^{-\frac{1}{2}}
    \|\tilde{u}(t-s+T_{m})\|_{L^\infty(\mathbb{R}^N)}ds\\
    &\leq e^{-\sigma(\gamma-1)(t-T_{m})}(\|\nabla v_0\|_{L^\infty(\mathbb{R}^N)}+\mu \bar{C}C_NC_1)
    +\mu \tilde{C}C_Ne^{-\sigma(\gamma-1) t}\int_0^{t-T_{m}} e^{-(\lambda-\sigma(\gamma-1))s}s^{-\frac{1}{2}}ds\\
    &\leq e^{-\sigma(\gamma-1)(t-T_{m})}(\|\nabla v_0\|_{L^\infty(\mathbb{R}^N)}+\mu \bar{C}C_NC_1+\mu \tilde{C}C_NC_3),
\end{align*}
where $C_N$ is defined in Lemma \ref{le2.1}, $C_1=\int_0^{+\infty}t^{-\frac{1}{2}}e^{-\lambda t}dt<+\infty$ and
$C_3=\int_0^{+\infty}e^{-(\lambda-\sigma(\gamma-1))s}s^{-\frac{1}{2}}ds<+\infty.$ Thus we have
\begin{align*}
    \|\nabla\tilde{v}(t)\|_{L^\infty(\mathbb{R}^N)}
    &\leq  C_{\tilde{v}}e^{-\sigma(\gamma-1)(t-T_{m})}
\end{align*}
for all $t\in(T_{m},\tilde{T})$ with $C_{\tilde{v}}=\|\nabla v_0\|_{L^\infty(\mathbb{R}^N)}
+\mu \bar{C}C_NC_1+\mu \tilde{C}C_NC_3$.

{\bf Claim 2.}
\textit{There exist} $\sigma\in(0,\min\{\frac{a}{2},
\frac{\lambda}{\gamma-1}\})$,
$\varepsilon\in(0,1)$ \textit{and} $\chi_*>0$ \textit{such that for} $\chi<\chi_*$ \textit{we obtain that}
\begin{align*}
    \|\tilde{u}(t)\|_{L^\infty(\mathbb{R}^N)}\leq\varepsilon \tilde{C} e^{-\sigma(\gamma-1)(t-T_{m})},
\end{align*}
\textit{for all} $t\in(T_{m},\tilde{T})$, \textit{where} $\tilde{C}=\frac{\varepsilon
a^{\frac{2-\gamma}{\gamma-1}}(a-\sigma)}{2b^{\frac{1}{\gamma-1}}}$.

According to \eqref{6.2} and using Lemma \ref{le2.1} and Lemma \ref{le2.2}, we deduce
\begin{align*}
    \|\tilde{u}(t)\|_{L^\infty(\mathbb{R}^N)}
    &\leq e^{-a(\gamma-1)(t-T_{m})}\|\tilde{u}(T_{m})\|_{L^\infty(\mathbb{R}^N)}\\
    &~~~+a(\gamma-1)\left(\frac{b}{a}\right)^{\frac{1}{\gamma-1}}
    \int_{T_{m}}^te^{-a(\gamma-1)(s-T_{m})}\|\tilde{u}^2(t-s+T_{m})\|_{L^\infty(\mathbb{R}^N)}ds\\
    &~~~+\frac{\chi N}{\sqrt{\pi}}\left(\frac{a}{b}\right)^{\frac{1}{\gamma-1}}
    \int_{T_{m}}^te^{-a(\gamma-1)(s-T_{m})}(s-T_{m})^{-\frac{1}{2}}\|\nabla\tilde{v}(t-s+T_{m})
    \|_{L^\infty(\mathbb{R}^N)}ds\\
    &~~~+\frac{\chi N}{\sqrt{\pi}}\int_{T_{m}}^te^{-a(\gamma-1)(s-T_{m})}(s-T_{m})^{-\frac{1}{2}}
    \|\tilde{u}(t-s+T_{m})\|_{L^\infty(\mathbb{R}^N)}
    \|\nabla\tilde{v}(t-s+T_{m})\|_{L^\infty(\mathbb{R}^N)}ds\\
    &\leq e^{-a(\gamma-1)(t-T_{m})}\frac{\varepsilon^2a^{\frac{2-\gamma}{\gamma-1}}(a-2\sigma)}
    {4b^{\frac{1}{\gamma-1}}} \\
    &~~~+ a(\gamma-1)\left(\frac{b}{a}\right)^{\frac{1}{\gamma-1}}{\tilde{C}}^2
    \int_{T_{m}}^te^{-a(\gamma-1)(s-T_{m})}e^{-2\sigma(\gamma-1)(t-s+T_{m})}ds\\
    &~~~+\frac{\chi N\tilde{C}C_{\tilde{v}}}{\sqrt{\pi}}\int_{T_{m}}^te^{-a(\gamma-1)(s-T_{m})}
    (s-T_{m})^{-\frac{1}{2}}e^{-2\sigma(\gamma-1)(t-s+T_{m})}ds\\
    &~~~+\left(\frac{a}{b}\right)^{\frac{1}{\gamma-1}}\frac{\chi NC_{\tilde{v}}}
    {\sqrt{\pi}}\int_{T_{m}}^te^{-a(\gamma-1)(s-T_{m})}(s-T_{m})^{-\frac{1}{2}}e^{-\sigma(\gamma-1)(t-s+T_{m})}ds\\
    &\leq e^{-a(\gamma-1)(t-T_{m})}\frac{\varepsilon^2a^{\frac{2-\gamma}{\gamma-1}}(a-2\sigma)}
    {4b^{\frac{1}{\gamma-1}}}+\frac{a\left(\frac{b}{a}\right)^{\frac{1}{\gamma-1}}{\tilde{C}}^2}
    {a-\sigma}e^{-\sigma(\gamma-1)(t-T_{m})}\\
    &~~~+\frac{\chi NC_{\tilde{v}}}{\sqrt{\pi}}
    (\tilde{C}+\left(\frac{a}{b}\right)^{\frac{1}{\gamma-1}})e^{-\sigma(\gamma-1) (t-T_{m})}\int_{T_{m}}^{t-T_{m}}
    e^{-(a(\gamma-1)-\sigma(\gamma-1))(s-T_{m})}(s-T_{m})^{-\frac{1}{2}}ds.
\end{align*}
Due to $a-\sigma>a-2\sigma>0$, we obtain
\begin{align*}
    \|\tilde{u}(t)\|_{L^\infty(\mathbb{R}^N)}
    \leq e^{-\sigma(\gamma-1)(t-T_{m})}
    \left[\frac{\varepsilon^2a^{\frac{2-\gamma}{\gamma-1}}(a-2\sigma)}
    {4b^{\frac{1}{\gamma-1}}}+\frac{a\left(\frac{b}{a}\right)^{\frac{1}{\gamma-1}}{\tilde{C}}^2}
    {a-\sigma}+\frac{\chi NC_{\tilde{v}}C_4}{\sqrt{\pi}}(\tilde{C}+\left(\frac{a}{b}\right)^{\frac{1}{\gamma-1}})\right],
\end{align*}
where $C_4=\int_0^{+\infty} e^{-(a(\gamma-1)-\sigma(\gamma-1))s}s^{-\frac{1}{2}}ds<+\infty$.

Noticing that $\tilde{C}=\frac{\varepsilon a^{\frac{2-\gamma}{\gamma-1}}(a-\sigma)}{2b^{\frac{1}{\gamma-1}}}$,
which implies
\begin{align*}
    \varepsilon \tilde{C}-\frac{\varepsilon^2a^{\frac{2-\gamma}{\gamma-1}}(a-2\sigma)}
    {4b^{\frac{1}{\gamma-1}}}-\frac{a\left(\frac{b}{a}\right)^{\frac{1}{\gamma-1}}{\tilde{C}}^2}
    {a-\sigma}=\frac{\varepsilon^2a^{\frac{2-\gamma}{\gamma-1}}\sigma}{4b^{\frac{1}{\gamma-1}}}>0.
\end{align*}
Then we can conclude that for all $t\in(T_{m},\tilde{T})$
\begin{align*}
    \|\tilde{u}(t)\|_{L^\infty(\mathbb{R}^N)}\leq\varepsilon \tilde{C} e^{-\sigma(\gamma-1)(t-T_{m})}
\end{align*}
for
\begin{align*}
    \chi<\chi_*=\min\{\chi_0,\chi_1,\frac{\varepsilon^2a^{\frac{2-\gamma}{\gamma-1}}\sigma\sqrt{\pi}}
    {4b^{\frac{1}{\gamma-1}}NC_{\tilde{v}}C_4(\tilde{C}+\left(\frac{a}{b}\right)^{\frac{1}{\gamma-1}})}\},
\end{align*}
where $\chi_0, \chi_1$ defined in \eqref{3.7} and \eqref{6.5}.
$\hfill\Box$

Now, we show that $\tilde{T}=+\infty$.
Claim 2 allows us to see
$\|\tilde{u}(t)\|_{L^{\infty}\left(\mathbb{R}^{N}\right)}\leq
\varepsilon \tilde{C}e^{-\sigma(\gamma-1)\left(t-T_{m}\right)}$
with $\varepsilon \in(0,1)$ for all $t \in\left(T_{m}, \tilde{T}\right)$.
Suppose $\tilde{T}<+\infty$, it contradicts the definition of $\tilde{T}$.
Hence $\tilde{T}=+\infty$. Then the large time behavior of $v(t)$ can be also obtained by
Claim 1 for all $t>T_{m}$.
The proof of Theorem \ref{th3.4} is finished.
$\hfill\Box$

\section*{Conflicts of Interest}
Authors have no conflict of interest to declare.



{\small

{Qingchun Li, Haomeng Chen}\\
{School of Mathematical Sciences,}\\
{Dalian University of Technology,}\\
{Dalian 116024, P. R. China}
\\
\\

\end{document}